\documentclass[11pt, reqno]{article}
\usepackage{amsfonts}
 \usepackage{pifont}
\usepackage{amsmath,amsfonts,amssymb,theorem,graphicx,listings,txfonts}
\usepackage{graphics}
 \usepackage{caption2}
\usepackage{flafter}
 \usepackage{indentfirst}
\usepackage{epsfig}
  \usepackage{mathrsfs}
\usepackage{subfigure}
 \usepackage{cite}
\usepackage{amsmath}
\usepackage{hyperref}
\usepackage{color}

\newtheorem{theorem}{Theorem}[section]
\newtheorem{lemma}[theorem]{Lemma}

\newcommand{\N}{\mbox{$\mathbb{N}$}}

\newcommand{\C}{\mbox{$\mathbb{C}$}}

\newcommand{\D}{\mbox{$\mathbb{D}$}}
\numberwithin{equation}{section}
\newcommand{\abs}[1]{\left\vert#1\right\vert}

\newtheorem{con}{Conjecture}[section]
\newtheorem{problem}{Problem}[section]
\newtheorem{dfn}[theorem]{Definition}

\newtheorem{rem}[theorem]{Remark}
\theoremstyle{example}

\theoremstyle{procedure}

\title{\bf{On a problem of Bharanedhar and Ponnusamy involving\vskip.05in planar harmonic mappings}}

\author
{\bf{Zhi-Gang Wang, Zhi-Hong Liu, Antti Rasila and Yong Sun}}

\date{}
\begin{document}
\maketitle \baselineskip=14.30pt
\begin{center}
\begin{quote}
 {{\bf Abstract.}
In this paper, we give a negative answer to a problem presented by
Bharanedhar and Ponnusamy (Rocky Mountain J. Math. 44: 753--777,
2014) concerning univalency of a class of harmonic mappings. More
precisely, we show that for all values of the involved parameter,
this class contains a non-univalent function. Moreover, several results on a new subclass of close-to-convex
harmonic mappings, which is motivated by work of Ponnusamy and
Sairam Kaliraj (Mediterr. J. Math. 12: 647--665, 2015), are
obtained.}
\end{quote}
\end{center}
\renewcommand{\thesection}{\arabic{section}}

\maketitle

\renewcommand{\thefootnote}{\fnsymbol{footnote}}

\footnotetext{\hspace*{-5mm}
\begin{tabular}{@{}r@{}p{16.2cm}@{}}
& Date: \date{\today}.
\\ &\vskip.03in\textit{Key Words and Phrases.} Planar harmonic mapping; univalent harmonic mapping;
close-to-convex harmonic mapping.
\\& \vskip.03in\textit{2010 Mathematics Subject Classification.} Primary 58E20;\ Secondary 30C55.
\\ & \vskip.03in \baselineskip=14.30pt {Z.-G. Wang was supported by the \textit{Natural
Science Foundation of Hunan Province} under Grant no. 2016JJ2036 and the \textit{National
Natural Science Foundation} under Grant no. 11301008. Z.-H. Liu was supported by the \textit{Foundation
of Educational Committee of Yunnan Province}
under Grant no. 2015Y456 and the \textit{Young and Middle-aged Academic Training Object
of Honghe University} under Grant no. 2014GG0102.
A. Rasila was partially supported by \textit{Academy of Finland} under Grant no. 289576.}
\end{tabular}}

\vskip.05in

\section{Introduction}

In this paper, we consider univalency criteria for complex-valued harmonic functions $f$ in the open unit disk $\D$. It is well-known that such functions can be written
as  $f=h+\overline{g}$, where $h$ and $g$ are analytic functions in $\D$. We call $h$ the analytic part and $g$ the co-analytic part of $f$, respectively. Let $\mathcal{H}$ be
the class of harmonic functions normalized by the conditions $f (0) = f_z(0)-1=0$, which have the
form
\begin{equation}\label{111}f (z)=z+\sum_{k=2}^{\infty}a_kz^k+\overline{\sum_{k=1}^{\infty}b_kz^k} \quad (z\in\D).\end{equation} Since the Jacobian of $f$ is
given by $\abs{h'}^2-\abs{g'}^2$, by Lewy's theorem (see
\cite{Lewy}), it is locally univalent and sense-preserving if and
only if $|g^{\prime}|<|h^{\prime}|$, or equivalently, the dilatation
$\omega={g^{\prime}}/{h^{\prime}}$ with
$h^{\prime}(z)\neq 0$ has the property $|\omega|<1$ in
$\D$. The subclass of
$\mathcal{H}$ that are harmonic, univalent and sense-preserving
in $\D$ is denoted by $\mathcal{S_{\mathcal{H}}}$. Univalent harmonic functions are also called harmonic mappings.

The classical family $\mathcal{S}$ of analytic univalent and
normalized functions in $\D$ is a subclass of
$\mathcal{S_\mathcal{H}}$ with $g(z)\equiv 0$. The family of all
functions $f\in\mathcal{S}_{\mathcal{H}}$ with the additional
property that $f_{\overline{z}}(0)=0$ is denoted by
$\mathcal{S}^{0}_{\mathcal{H}}$. There exist reciprocal transformations
 between the classes $\mathcal{S}_{\mathcal{H}}$
and $\mathcal{S}^{0}_{\mathcal{H}}$ (see \cite{cs,d}). Observe that the family $\mathcal{S}^0_{\mathcal{H}}$
is compact and normal, but the family $\mathcal{S}_{\mathcal{H}}$
is not compact.  For recent results involving univalent harmonic mappings, we refer to \cite{bjj,bl,bp,cprw,k11,kpm,m5,nr,ps,ps2,ps1,sjr,srj,wll,wsj1}, and the references therein.

A domain $\Omega$ is said to be close-to-convex if $\C\backslash\Omega$ can be represented as a union of non-intersecting
half-lines. Following the result due to Kaplan \cite{k}, an analytic function $F$ is called close-to-convex if there exists a
univalent convex analytic function $\phi$ defined in $\D$ such that
$${\rm Re}\left(\frac{F'(z)}{\phi'(z)}\right)>0\quad (z\in\D).$$ Furthermore, a planar harmonic
mapping $f:\D\rightarrow\C$ is close-to-convex if it is injective
and $f(\D)$ is a close-to-convex domain, we denote by $\mathcal{C}_\mathcal{H}^0$ the class of close-to-convex harmonic mappings.

This paper is organized as follows. In Section 2, we give a negative
answer to a problem posed by Bharanedhar and Ponnusamy in \cite{bp}.  In Section 3, we study a subclass of
close-to-convex harmonic mappings, which is motivated by work of
Ponnusamy and Sairam Kaliraj \cite{ps1}. Coefficient estimates, a
growth theorem, a covering theorem and an area theorem, for mappings of this class,
are obtained.

\vskip.08in

\section{A problem of Bharanedhar and Ponnusamy}

Recently, Mocanu \cite{m5} proposed the following conjecture involving the
univalency of planar harmonic mappings.

\begin{con}\label{c111} Let $$\mathcal{M}=\left\{f=h+\overline{g}\in \mathcal{H}:\ g'=zh'\ {\it and}\
{\rm Re}\left(1+\frac{zh''(z)}{h'(z)}\right)>-\frac{1}{2}\quad(z\in\D)\right\}.$$ Then $\mathcal{M}\subset\mathcal{S}_\mathcal{H}^0$.
\end{con}

By applying the close-to-convexity criterion for analytic functions due to Kaplan \cite{k}, Bshouty and Lyzzaik \cite{bl} have solved the above conjecture by establishing
the following stronger result:
\vskip.10in
\noindent{\bf Theorem A}\ \ $\mathcal{M}\subset\mathcal{C}_\mathcal{H}^0$.
\vskip.10in

Later, Ponnusamy and Sairam Kaliraj \cite[Theorem 4.1]{ps1} generalized Theorem A, under the assumption that the analytic dilatation $\omega$ satisfies the condition
$${\rm Re}\left(\frac{\lambda z \omega'(z)}{1-\lambda\omega(z)}\right)>-\frac{1}{2}$$ for all $\lambda$ such that $\abs{\lambda}=1$.
In particular, for $$\omega(z)=\lambda kz^n \quad \left(\abs{\lambda}=1;\,0<k\leq \frac{1}{2n-1};\,n\in\N:=\{1,2,3,\ldots\}\right),$$ they gave the following result:

\vskip.10in
\noindent{\bf Theorem B}\ \ \textit{Suppose that $h$ and $g$ are analytic in $\D$ such that $${\rm Re}\left(1+\frac{zh''(z)}{h'(z)}\right)>-\frac{1}{2},$$ and $$
g'(z)=\lambda kz^nh'(z)
\quad \left(n\in\N;\, \abs{\lambda}=1;\, 0<k\leq \frac{1}{2n-1}\right).$$ Then $f=h+\overline{g}$ is univalent and close-to-convex in $\D$.}

Motivated by Theorem B, we introduce the following natural class of close-to-convex harmonic mappings, which will be studied in Section 3. Note that for $n=1$, we have
the class $\mathcal{M}(\alpha,\zeta)$, which was studied in \cite{sjr}.

\begin{dfn}
{\rm A harmonic mapping $f=h+\overline{g}\in\mathcal{H}$ is said to be in the class $\mathcal{M}(\alpha,\zeta,n)$ if
$h$ and $g$ satisfy the
conditions \begin{equation}\label{401}{\rm Re}\left(1+\frac{zh''(z)}{h'(z)}\right)>\alpha \quad\left(-\frac{1}{2}\leq\alpha<1\right),\end{equation} and
\begin{equation}\label{402}
g'(z)=\zeta z^nh'(z)
\quad \left(\zeta\in\C\ {\rm with}\ \abs{\zeta}\leq \frac{1}{2n-1};\, n\in\N\right).\end{equation}}
\end{dfn}

In 1995, Ponnusamy and Rajasekaran \cite{pr} derived the following starlikeness criterion for analytic functions.
\vskip.10in
\noindent{\bf Theorem C}\ \ \textit{Suppose that $F$ is a normalized analytic function in $\D$. If $F$ satisfies the condition $${\rm Re}\left(1+\frac{zF''(z)}{F'(z)}\right)<\beta\quad\left(1<\beta\leq \frac{3}{2}\right),$$ then $F$ is univalent and starlike in $\D$,
i.e. $F(\D)$ is a domain
starlike with respect to the origin.}

Motivated essentially by Theorems A and C, Bharanedhar and Ponnusamy \cite[Problem 1, p. 763]{bp}
posed the following problem, which we present here in a slightly modified form:

\begin{problem}\label{c222} For $\beta\in(1,3/2)$, define $$\mathcal{P}(\beta)=\left\{f=h+\overline{g}\in \mathcal{H}:\ g'=zh'\ {\it and}\
{\rm Re}\left(1+\frac{zh''(z)}{h'(z)}\right)<\beta\quad(z\in\D)\right\}.$$ Determine $\inf\left\{\beta\in(1,3/2):\, \mathcal{P}(\beta)\subset \mathcal{S}^0_{\mathcal{H}}\right\}$.
\end{problem}

Let us recall the following result of Bshouty and Lyzzaik \cite{bl}:

\vskip.10in
\noindent{\bf Theorem D}\ \ \textit{Suppose that $0\leq\lambda<1/2$. Let $f=h+\overline{g}$ be the harmonic polynomial mapping with $$h(z)=z-\lambda z^2 \ \ {\it {and}}\ \ g(z)=\frac
{z^2}{2}-\frac{2\lambda z^3}{3}.$$
If\ \ $0\leq\lambda\leq 3/10$, then $f$ is univalent in $\D$. But for $3/10<\lambda<1/2$, $f$ is not univalent in $\D$.}

\begin{rem}{\rm In view of Theorem D, we see that $\beta$ can be restricted to the value on the interval $(1,11/8]$, since $$\sup\limits_{z\in{\scriptsize \D}}\left\{{\rm Re}\left(1+\frac{zh''(z)}{h'(z)}\right)\right\}=\frac{11}{8}$$
for} $$h(z)=z-\frac{3}{10}z^2.$$
\end{rem}

Now, we are ready to give a counterexample, which shows that for all $\beta\in(1,11/8]$, the class $\mathcal{P}(\beta)$ of Problem \ref{c222} contains a non-univalent function.

Consider the harmonic function given by $f_{\gamma}=h+\overline{g}\in\mathcal{H}$, where
$$h(z)=\frac{1}{\gamma}\left[1-(1-z)^{\gamma}\right]\quad\left(1<\gamma\leq\frac{7}{4}\right),$$ and

$$g(z)=\frac{1}{\gamma(1+\gamma)}\left[1-(1+\gamma z)(1-z)^{\gamma}\right]\quad\left(1<\gamma\leq\frac{7}{4}\right).$$
Clearly, we have $g'=zh'$. It follows that $$1+\frac{zh''(z)}{h'(z)}=\frac{1-\gamma z}{1-z},$$ and therefore  $${\rm Re}\left(1+\frac{zh''(z)}{h'(z)}\right)<\frac{1+\gamma}{2}\quad\left(1<\frac{1+\gamma}{2}\leq\frac{11}{8}\right).$$ That is,
 $$f_{\gamma}=h+\overline{g}\in\mathcal{P}((1+\gamma)/2)\subset\mathcal{P}(\beta).$$

In what follows, we shall prove that the function $f_{\gamma}$ is not univalent in $\D$. It is easy to verify that both the analytic and co-analytic parts of $f_{\gamma}$ have real
coefficients, and thus, $f_{\gamma}(z)=\overline{f_{\gamma}(\overline{z})}$ for all $z\in\D$. In particular,
$${\rm Re}\left(f_{\gamma}\left(r e^{i \theta}\right)\right)={\rm Re}\left(f_{\gamma}\left(r e^{-i \theta}\right)\right)
$$
for some $r\in (0,1)$ and $\theta\in (-\pi,0)\cup(0,\pi)$.
It suffices to show that there exist $r_{0}\in (0,1)$ and $\theta_{0}\in (-\pi,0)\cup(0,\pi)$ such that
$${\rm Im} \left(f_{\gamma}\left(r_{0} e^{i \theta_{0}}\right)\right)={\rm Im}\left(f_{\gamma}\left(r_{0} e^{-i \theta_{0}}\right)\right)=0.
$$
In view of the relation
\begin{equation*}
{\rm Im}\left( f_{\gamma}(z)\right)={\rm Im}\left(h(z)-g(z)\right)
={\rm Im}\left(\frac{1-(1-z)^{\gamma+1}}{\gamma+1}\right)
=-{\rm Im}\left(\frac{e^{(\gamma+1)\log(1-z)}}{\gamma+1}\right),
\end{equation*}
we see that
\begin{equation*}
\begin{split}
{\rm Im} \left(f_{\gamma}\left(r e^{i\theta}\right)\right)&=-{\rm Im}\left(\frac{e^{(\gamma+1)\log\left(1-r e^{i\theta}\right)}}{\gamma+1}\right)\\
&=-\frac{e^{(\gamma+1)\log|1-r e^{i\theta}|}}{\gamma+1}\sin \left[(\gamma+1)\arg\left(1-r e^{i\theta}\right)\right],
\end{split}
\end{equation*}
and
\begin{equation*}
\begin{split}
-{\rm Im}\left(f_{\gamma}\left(r e^{-i\theta}\right)\right)=\frac{e^{(\gamma+1)\log|1-r e^{-i\theta}|}}{\gamma+1}\sin \left[(\gamma+1)\arg\left(1-r e^{-i\theta}\right)\right]={\rm Im}\left(f_{\gamma}\left(r e^{i\theta}\right)\right).
\end{split}
\end{equation*}
By noting that $$\arg\left(1-r e^{i\theta}\right)\in\left(-\frac{\pi}{2},0\right)\cup \left(0,\frac{\pi}{2}\right),
$$ we deduce that for each $1<\gamma\leq7/4$, there exist $r_{0}\in (0,1)$ and $\theta_{0}\in (-\pi,0)\cup(0,\pi)$ such that
$$\sin\left[(\gamma+1)\arg\left(1-r_{0} e^{i\theta_{0}}\right)\right]=0.
$$
It follows that $${\rm Im}\left(f_{\gamma}\left(r_{0} e^{i \theta_{0}}\right)\right)={\rm Im} \left(f_{\gamma}\left(r_{0} e^{-i \theta_{0}}\right)\right)=0.$$ Therefore,  there exist two distinct points $z_1 = r_{0} e^{i \theta_{0}}$ and $z_2 = r_{0} e^{-i \theta_{0}}$ in $\D$ such that $f_{\gamma}(z_1) = f_{\gamma}(z_2)$, which shows that the function $f_{\gamma}(z)$ is not univalent in $\D$.
Thus, we conclude that the conditions given in Problem \ref{c222} are not satisfied for any $\beta\in(1,11/8]$.

The image domain of $f_{\gamma}$ for $\gamma=5/4$ is given in Figures 1 and 2 to illustrate our counterexample.
\begin{figure}[htbp]
\begin{minipage}{0.45\linewidth}\label{f111}
    \centering
\includegraphics[width=2.4in]{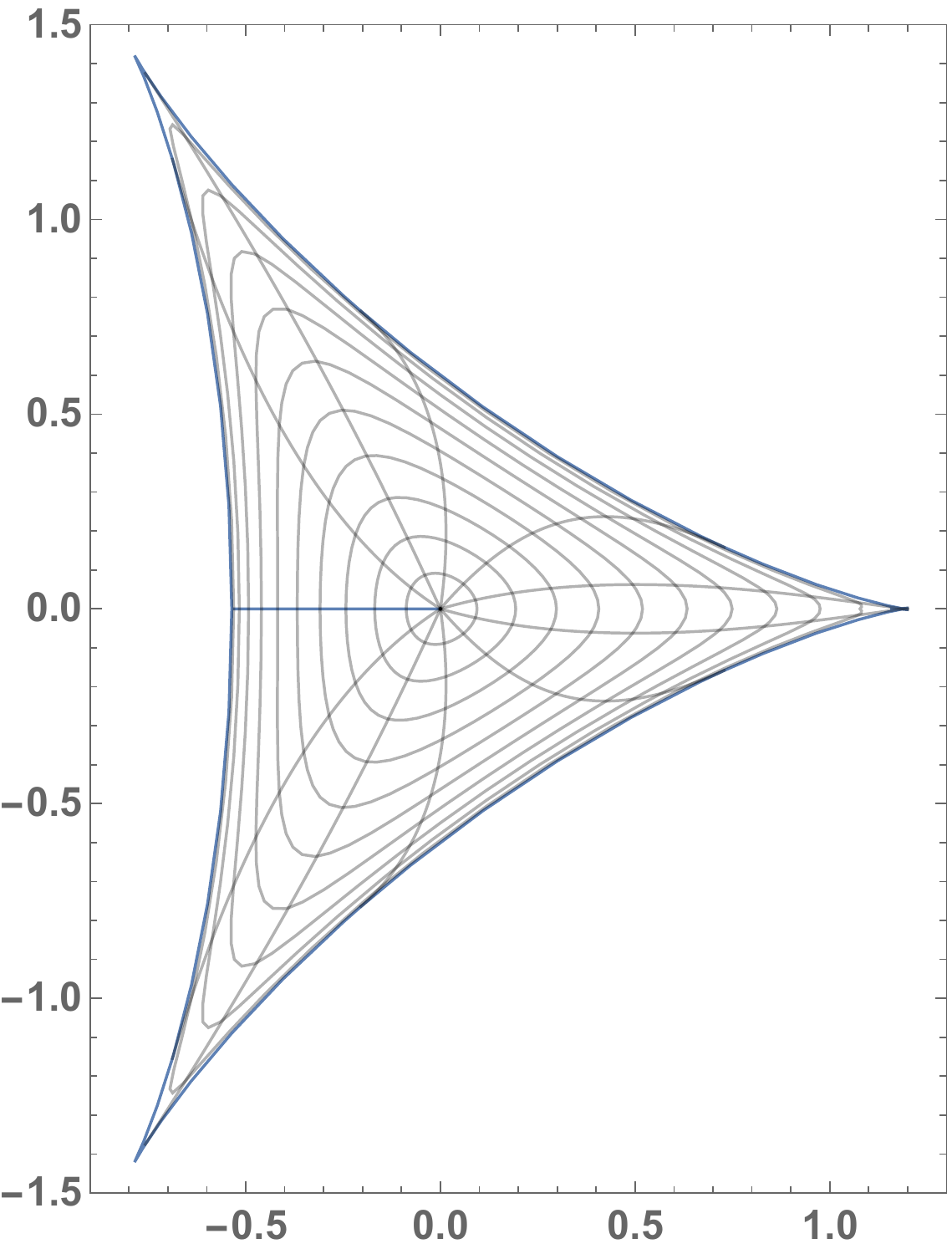}
\caption{The image of the mapping $f_{5/4}$.}\end{minipage}
\hspace{-1ex}
\begin{minipage}{0.45\linewidth}\label{f22}
    \centering
\includegraphics[width=3.3in]{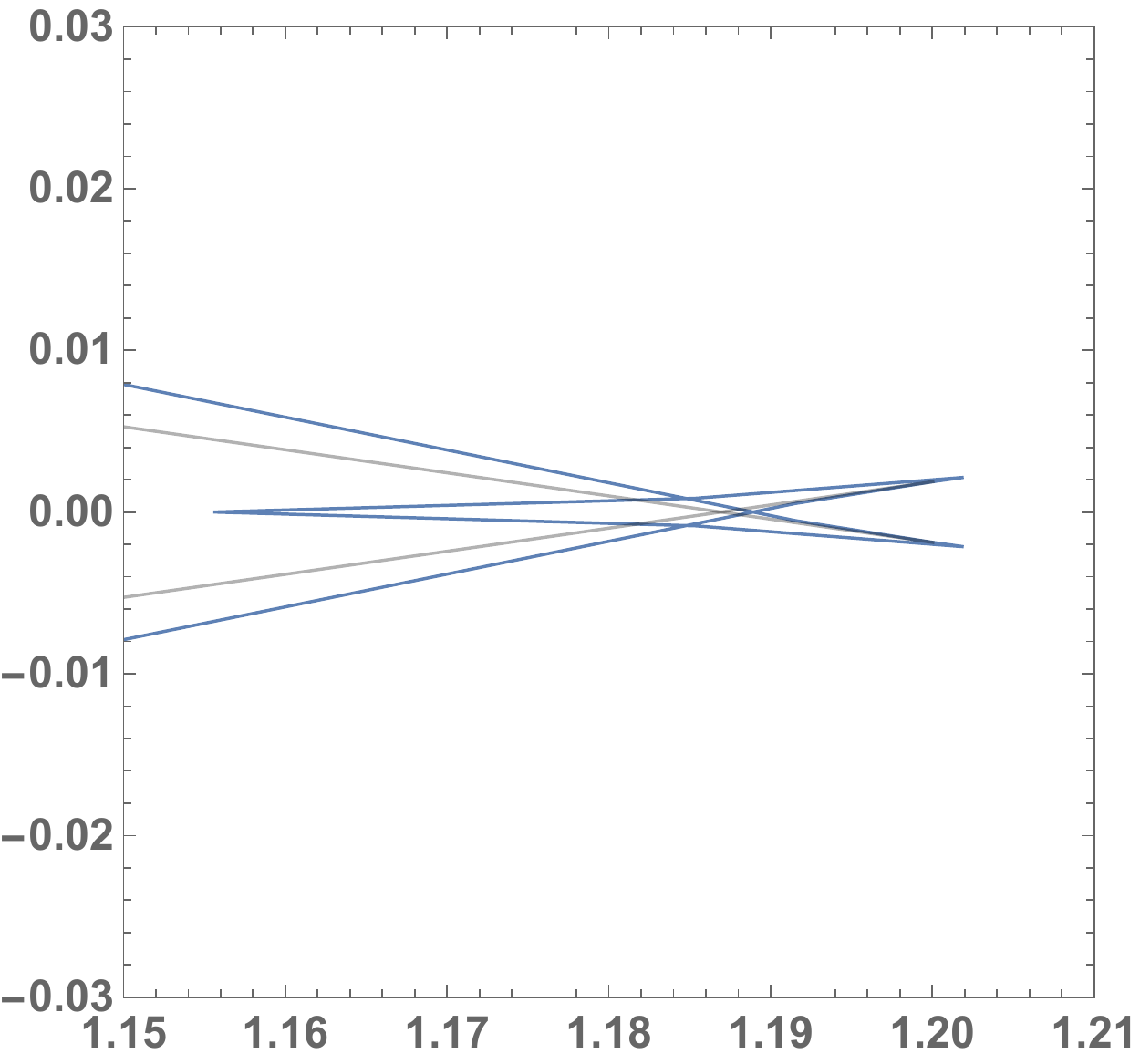}
\caption{An enlarged view of right cusp of image of $f_{5/4}$.}
\end{minipage}
\end{figure}

\vskip.08in

\section{The subclass $\mathcal{M}(\alpha,\zeta,n)$ of close-to-convex harmonic mappings}

Recall the following lemma, due to Suffridge \cite{s2005}, which will be required in the proof of Theorem \ref{t41}.

\begin{lemma}\label{l41}
If $h(z)=z+\sum_{k=2}^{\infty}a_kz^k$ satisfies the condition \eqref{401}, then
\begin{equation}\label{41}\abs{a_k}\leq\frac{1}{k!}\prod_{j=2}^{k}(j-2\alpha)\quad(k\in\N\setminus\{1\}),\end{equation}
with the extremal function given by
$$h(z)=\int_0^z\frac{dt}{(1-\delta t)^{2-2\alpha}}\quad(\abs{\delta}=1;\, z\in\D).$$
\end{lemma}

We now derive the coefficient estimates for the class $\mathcal{M}(\alpha,\zeta,n)$.

\begin{theorem}\label{t41}
Let $f=h+\overline{g}\in\mathcal{M}(\alpha,\zeta,n)$ be of the form \eqref{111}. Then
the coefficients $a_k\ (k\in\N\setminus\{1\})$ of $h$ satisfy \eqref{41}.
Moreover, the coefficients $b_k\ (k=n+1, n+2, \cdots; n\in\N)$ of $g$ satisfy 
\[
\abs{b_{n+1}}\leq\frac{\abs{\zeta}}{n+1}\ \ \text{ and }\ \, \ \, \abs{b_{k+n}}\leq\frac{\abs{\zeta}}{(k+n)(k-1)!}\prod_{j=2}^{k}(j-2\alpha)\quad(k\in\N\setminus\{1\};\, n\in\N).
\]
The bounds are sharp for the extremal function given by
$$f(z)=\int_0^z\frac{dt}{(1-\delta t)^{2-2\alpha}}+\overline{\int_0^z\frac{\zeta t^n}{(1-\delta t)^{2-2\alpha}}dt}\quad(\abs{\delta}=1;\, z\in\D).$$
\end{theorem}

\noindent{\bf Proof.}
By equating the coefficients of $z^{k+n-1}$ in both sides of \eqref{402}, we see that \begin{equation}\label{45}(k+n)b_{k+n}=\zeta k a_k\quad (k,\,n\in\N;\, a_1=1).\end{equation} In view of Lemma \ref{l41} and \eqref{45}, we get the desired result of Theorem \ref{t41}.
\hfill\rule{1.5mm}{3mm}

\begin{theorem}\label{t401}
Let $f\in\mathcal{M}(\alpha,\zeta,n)$ with $0\leq\alpha<1$ and $0\leq\zeta<\frac{1}{2n-1}\ (n\in\N)$.
Then
\begin{equation}\label{411}
\Phi(r;\alpha,\zeta,n)\leq \abs{f(z)}\leq \Psi(r;\alpha,\zeta,n)\quad(r=|z|<1),
\end{equation}
where
\begin{equation*}
\Phi(r;\alpha,\zeta,n)
=\left\{\begin{array}{ll}
\log (1+r)-\displaystyle\frac{\zeta\,  r^{n+1} \, _2F_1(1,\, n+1;\, n+2;\, -r)}{n+1}  &(\alpha=1/2), \\ \\
\displaystyle\frac{(1+r)^{2 \alpha -1}-1}{2 \alpha -1}-\frac{\zeta\,  r^{n+1} \, _2F_1(n+1,\, 2-2 \alpha;\, n+2;\, -r)}{n+1}  &(\alpha\neq 1/2),
\end{array}\right.
\end{equation*}
and
\begin{equation*}
\Psi(r;\alpha,\zeta,n)
=\left\{\begin{array}{ll}-\log (1-r)+\displaystyle\frac{\zeta\,  r^{n+1} \, _2F_1(1,\, n+1;\, n+2;\, r)}{n+1}  &(\alpha=1/2), \\ \\
\displaystyle\frac{1-(1-r)^{2 \alpha -1}}{2 \alpha -1}+\frac{\zeta\,  r^{n+1} \, _2F_1(n+1,\, 2-2 \alpha;\, n+2;\, r)}{n+1}  &(\alpha\neq 1/2).
\end{array}\right.
\end{equation*}
All these bounds are sharp, the extremal function is $f_{\alpha,\zeta,n}=h_{\alpha}+\overline{g_{\alpha,\zeta,n}}$
or its rotations, where
\begin{equation}\label{412}
f_{\alpha,\zeta,n}(z)
=\left\{\begin{array}{ll}
-\log(1-z)+\displaystyle\overline{\frac{\zeta\,  z^{n+1} \, _2F_1(1,\, n+1;\, n+2;\, z)}{n+1}}  &(\alpha=1/2), \\ \\
\displaystyle\frac{1-(1-z)^{2\alpha-1}}{2\alpha-1}
+\overline{\frac{\zeta\,  z^{n+1} \, _2F_1(n+1,\, 2-2 \alpha;\, n+2;\, z)}{n+1}}  &(\alpha\neq 1/2).
\end{array}\right.
\end{equation}
\end{theorem}

\noindent{\bf Proof.}
Assume that $f=h+\overline{g}\in\mathcal{M}(\alpha,\zeta,n)$.
Also, let $\Gamma$ be the line segment joining $0$ and $z$. Then
\begin{equation}\label{413}
\abs{f(z)}=\abs{\int_{\Gamma}\frac{\partial f}{\partial \xi}d\xi
+\frac{\partial f}{\partial \overline{\xi}}d\overline{\xi}}
\leq \int_{\Gamma}\left(\abs{h'(\xi)}+\abs{g'(\xi)}\right)\abs{d\xi}
= \int_{\Gamma}\left(1+\abs{\zeta}\abs{\xi}^n\right)\abs{h'(\xi)}\abs{d\xi}.
\end{equation}
Moreover, let $\widetilde{\Gamma}$ be the preimage under $f$ of the line segment joining $0$ and $f(z)$.  Then we obtain
\begin{equation}\label{414}
\abs{f(z)}=\int_{\widetilde{\Gamma}}\abs{\frac{\partial f}{\partial \xi}d\xi
+\frac{\partial f}{\partial \overline{\xi}}d\overline{\xi}}
\geq \int_{\widetilde{\Gamma}}\left(\abs{h'(\xi)}-\abs{g'(\xi)}\right)\abs{d\xi}
= \int_{\widetilde{\Gamma}}\left(1-|\zeta||\xi|^n\right)\abs{h'(\xi)}\abs{d\xi}.
\end{equation}
By observing that $h$ is a convex analytic function of order $\alpha\ (0\leq\alpha<1)$,
it follows that
\begin{equation}\label{415}
\frac{1}{(1+r)^{2(1-\alpha)}} \leq |h'(z)| \leq \frac{1}{(1-r)^{2(1-\alpha)}}
\quad (|z|=r<1).
\end{equation}
By virtue of \eqref{413}, \eqref{414} and \eqref{415}, we see that
\begin{equation*}
\Phi(r;\alpha,\zeta,n):=\int_{0}^{r}\frac{(1-|\zeta|\rho^n)d\rho}{(1+\rho)^{2(1-\alpha)}}
\leq |f(z)| \leq
\int_{0}^{r}\frac{(1+|\zeta|\rho^n)d\rho}{(1-\rho)^{2(1-\alpha)}}=:\Psi(r;\alpha,\zeta,n),
\end{equation*}
which yields the desired inequalities \eqref{411}.

Now, we shall prove the sharpness of the result. We only need to show that $f_{\alpha,\zeta,n}$
defined by \eqref{412} belongs to the class $\mathcal{M}(\alpha,\zeta,n)$ for each $\alpha\in[0,1)$. Suppose that
\begin{equation*}
h_{\alpha}(z)
=\left\{\begin{array}{ll}-\log(1-z)  &(\alpha=1/2), \\ \\
\displaystyle\frac{1-(1-z)^{2\alpha-1}}{2\alpha-1}  &(\alpha\neq 1/2).
\end{array}\right.
\end{equation*}
Then, we find that $h_{\alpha}(z)$ satisfies the inequality \eqref{401} and the relation
$g'_{\alpha,\zeta, n}(z)=\zeta z^n h'_{\alpha}(z)$ for each $\alpha\in[0,1)$.
Moreover, for $0\leq\alpha<1$, $0<\zeta<1/(2n-1)$ with $n\in\N$,  and $0<r<1$, it is easy to see that
\begin{equation*}
f_{\alpha,\zeta,n}(-r)=-\Phi(r;\alpha,\zeta,n)
\quad {\rm and} \quad
f_{\alpha,\zeta,n}(r)=\Psi(r;\alpha,\zeta,n),
\end{equation*}
and therefore,
\begin{equation*}
\big|f_{\alpha,\zeta,n}(-r)\big|=\Phi(r;\alpha,\zeta,n)
\quad {\rm and} \quad
\big|f_{\alpha,\zeta,n}(r)\big|=\Psi(r;\alpha,\zeta,n).
\end{equation*}
This shows that the bounds are sharp.
\hfill\rule{1.5mm}{3mm}

Next, we consider a covering theorem for functions in the class $\mathcal{M}(\alpha,\zeta,n)$.

\begin{theorem}\label{t42}
Let $f\in\mathcal{M}(\alpha,\zeta,n)$ with $0\leq\alpha<1$ and $0\leq\zeta<\frac{1}{2n-1}\ (n\in\N)$.
Then the range $f(\D)$ contains the disk
\begin{equation*}
|\omega|<r(\alpha,\zeta,n)
=\left\{\begin{array}{ll}
\log 2-\displaystyle\frac{\zeta \, _2F_1(1,\, n+1;\, n+2;\, -1)}{n+1}  &(\alpha=1/2), \\ \\
\displaystyle\frac{2^{2 \alpha -1}-1}{2 \alpha -1}-\frac{\zeta \, _2F_1(n+1,\, 2-2 \alpha;\, n+2;\, -1)}{n+1}  &(\alpha\neq 1/2).
\end{array}\right.
\end{equation*}
The bounds are sharp for the function $f_{\alpha,\zeta,n}=h_{\alpha}+\overline{g_{\alpha,\zeta,n}}$
given by \eqref{412} or its rotations.
\end{theorem}

\noindent{\bf Proof.}
By putting $r\rightarrow1^{-}$ in the lower bound for $|f(z)|$ in Theorem \ref{t401}, we get the desired result.
The sharpness is similar to that of Theorem \ref{t401}, we choose to omit the details.
\hfill\rule{1.5mm}{3mm}

Now, we consider the area theorem of the mappings belonging to the
class $\mathcal{M}(\alpha,\zeta,n)$. Let us denote
$\mathcal{A}\left(f(\D_{r})\right)$ by the area of $f(\D_{r})$,
where $\D_{r}:=r\D$ for $0<r<1$.

\begin{theorem}\label{t43}
Let $f\in\mathcal{M}(\alpha,\zeta,n)$ with $0\leq\alpha<1$.
Then, for $0<r<1$, $\mathcal{A}\left(f(\D_{r})\right)$ satisfies the inequalities
\begin{equation}\label{421}
2\pi\int_{0}^{r}\frac{\rho(1-|\zeta|^{2}\rho^{2n})}{(1+\rho)^{4(1-\alpha)}}d\rho
\leq \mathcal{A}\left(f(\D_{r})\right)\leq
2\pi\int_{0}^{r}\frac{\rho(1-|\zeta|^{2}\rho^{2n})}{(1-\rho)^{4(1-\alpha)}}d\rho.
\end{equation}
\end{theorem}

\noindent{\bf Proof.}
Let $f=h+\overline{g}\in\mathcal{M}(\alpha,\zeta,n)$. Then for $0<r<1$, we see that
\begin{equation}\label{4012}
\mathcal{A}\left(f(\D_{r})\right)=\iint_{\D_{r}}\left(|h'(z)|^{2}-|g'(z)|^{2}\right)\,dx\,dy
=\iint_{\D_{r}}\left(1-|\zeta|^{2}|z|^{2n}\right)|h'(z)|^{2}\,dx\,dy.
\end{equation}
By observing that $h$ is a convex analytic function of order $\alpha\ (0\leq\alpha<1)$,
in view of \eqref{415} and \eqref{4012}, we obtain the desired inequalities \eqref{421} of Theorem \ref{t43}.
\hfill\rule{1.5mm}{3mm}

\begin{rem}
{\rm By setting $n=1$ in Theorems \ref{t41}, \ref{t401}, \ref{t42} and \ref{t43}, respectively, we get the corresponding results obtained in \cite{sjr}.}
\end{rem}

Finally, we discuss the radius of close-to-convexity of a certain class harmonic mappings related to the class $\mathcal{M}(\alpha,\zeta,n)$.
The following lemma due to Clunie and Sheil-Small \cite{cs} will be required in the proof of Theorem \ref{t46}.

\begin{lemma}\label{lemma1} 
If $h,\,g$ are analytic in $\D$ with $\abs{h'(0)}>\abs{g'(0)}$, and $h+\lambda g$ is close-to-convex for each $\lambda\,
(\abs{\lambda}=1)$, then $f=h+\overline{g}$ is harmonic close-to-convex in $\D$.
\end{lemma}

\begin{theorem}\label{t46}
Suppose that $f=h+\overline{g}$ satisfies the inequality \eqref{401} with $-1/2<\alpha<0$. If $g'(z)=z^nh'(z)$ with $n\in\N\setminus\{1\}$,
then $f$ is close-to-convex in the disk
\begin{equation*}
\abs z<\sqrt[n]{\frac{1+2\alpha}{1+2n+2\alpha}}\quad(n\in\N\setminus\{1\}).
\end{equation*}
\end{theorem}

\noindent{\bf Proof.}
Suppose that $F_{\lambda}(z)=h(z)-\lambda g(z)$ with $\abs{\lambda}=1$. It follows that
$${\rm Re}\left(1+\frac{zF_{\lambda}''(z)}{F_{\lambda}'(z)}\right)={\rm Re}\left(1+\frac{zh''(z)}{h'(z)}\right)+n{\rm Re}\left(\frac{\lambda z^n}{\lambda z^n-1}\right)={\rm Re}\left(1+\frac{zh''(z)}{h'(z)}\right)+\frac{n}{2}\left(1-\frac{1-\abs{\lambda z^n}^2}{\left(1-\lambda z^n\right)\left(1-\overline{\lambda z^n}\right)}\right).$$
For $z=re^{i\theta}\ (0<r<1)$, we see that
$$\frac{n}{2}\left(1-\frac{1-\abs{\lambda z^n}^2}{\left(1-\lambda z^n\right)\left(1-\overline{\lambda z^n}\right)}\right)=\frac{n}{2}\left(1-\frac{1-r^{2n}}{1+r^{2n}-2{\rm Re}(\lambda z^{n})}\right)\geq-\frac{nr^n}{1-r^n}.$$
Thus,
$$\int_{\theta_1}^{\theta_2}{\rm Re}\left(1+\frac{zF_{\lambda}''(z)}{F_{\lambda}'(z)}\right)d\theta>\int_{\theta_1}^{\theta_2}\left(\alpha
-\frac{nr^n}{1-r^n}\right)d\theta=\left(\alpha
-\frac{nr^n}{1-r^n}\right)\left(\theta_2-\theta_1\right)>-\pi\quad(\theta_1<\theta_2<\theta_1+2\pi)$$ for $$\abs z=r<\sqrt[n]{\frac{1+2\alpha}{1+2n+2\alpha}}=:r(\alpha,n).$$
By Lemma \ref{lemma1} and Kaplan's close-to-convexity criterion for analytic functions (see \cite{k}), we deduce that $f$ is close-to-convex in the disk $\abs z<r(\alpha,n)$.
\hfill\rule{1.5mm}{3mm}

 \vskip.20in

\noindent\textbf{\Large Acknowledgments}

\vskip.10in

The authors would like to thank the referees and Prof. S. Ponnusamy for their valuable
comments and suggestions, which essentially improved the quality of this paper.

\vskip.18in

\textbf{Zhi-Gang Wang}
\vskip.05in {School of Mathematics and Computing Science, Hunan
First Normal University, Changsha 410205, Hunan, P.
R. China.}
\vskip.03in
\textit{E-mail address}: wangmath$@$163.com

\vskip.10in
\textbf{Zhi-Hong Liu}
\vskip.05in {School of Mathematics and Econometrics, Hunan University, Changsha 410082, Hunan, P.
R. China.}
\vskip.03in
\textit{E-mail address}: liuzhihongmath$@$163.com

\vskip.10in
\textbf{Antti Rasila} \vskip.05in{Department of
Mathematics and Systems Analysis, Aalto University, P. O. Box 11100,
FI-00076 Aalto, Finland.} \vskip.03in \textit{E-mail address}:
antti.rasila$@$iki.fi

\vskip.10in
\textbf{Yong Sun}
\vskip.05in{School of Science, Hunan Institute of Engineering, Xiangtan 411104, Hunan, P.
R. China.}
\vskip.03in
\textit{E-mail address}: yongsun2008$@$foxmail.com

\begin{thebibliography}{99}


\bibitem{bjj}
D. Bshouty, S. S. Joshi and S. B. Joshi, On close-to-convex harmonic mappings,  \textit{Complex Var. Elliptic Equ.} {\bf 58} (2013), 1195--1199.




\bibitem{bl}
D. Bshouty and A. Lyzzaik, Close-to-convexity criteria for planar
harmonic mappings, {\it Complex Anal. Oper. Theory} {\bf 5} (2011),
767--774.





\bibitem{bp}
S. V. Bharanedhar and S. Ponnusamy, Coefficient conditions for harmonic univalent mappings and hypergeometric mappings, \textit{Rocky Mountain J. Math.} {\bf 44} (2014), 753--777.





\bibitem{cprw}
 S. Chen, S. Ponnusamy,  A. Rasila and X. Wang,  Linear connectivity, Schwarz-Pick lemma and univalency criteria for planar harmonic mapping, \textit{Acta Math. Sin. (Engl. Ser.)} {\bf 32} (2016), 297--308.


\bibitem{cs}
J. Clunie and T. Sheil-Small, Harmonic univalent functions, \textit{
Ann. Acad. Sci. Fenn. Ser. A. I Math.} {\bf 9} (1984), 3--25.





\bibitem{d}
P. Duren, {\it Harmonic mappings in the plane}, Cambridge University
Press, Cambridge, 2004.



\bibitem{k11}
D. Kalaj, Quasiconformal harmonic mappings and close-to-convex domains, \textit{Filomat} \textbf{24} (2010), 63--68.





\bibitem{kpm}
D. Kalaj, S. Ponnusamy and M. Vuorinen, Radius of close-to-convexity and fully starlikeness of harmonic mappings, \textit{Complex Var. Elliptic Equ.} {\bf 59} (2014), 539--552.


\bibitem{k}
W. Kaplan, Close-to-convex schlicht functions, {\it Mich. Math. J.}
{\bf 1} (1952), 169--185.



\bibitem{Lewy}
H. Lewy, On the non-vanishing of the Jacobian in certain one-to-one
mappings, \textit{Bull. Amer. Math. Soc.} \textbf{42} (1936),
689--692.






\bibitem{m5}
P. T. Mocanu, Injectivity conditions in the complex plane, {\it
Complex Anal. Oper. Theory} {\bf 5} (2011), 759--766.



\bibitem{nr}
S. Nagpal and V. Ravichandran, Starlikeness, convexity and close-to-convexity of harmonic mappings,
\textit{Current topics in pure and computational complex analysis}, 201--214, Trends Math., Birkh\"{a}user/Springer, New Delhi, 2014.




\bibitem{pr}
S. Ponnusamy and S. Rajasekaran, New sufficient conditions for starlike and univalent functions, \textit{Soochow J. Math.} {\bf 21} (1995), 193--201.


\bibitem{ps}
S. Ponnusamy and A. Sairam Kaliraj, On harmonic close-to-convex functions, \textit{Comput. Methods Funct. Theory} {\bf 12} (2012), 669--685.


\bibitem{ps2}
S. Ponnusamy and A. Sairam Kaliraj, Univalent harmonic mappings convex in one direction, \textit{Anal. Math. Phys.} {\bf 4} (2014), 221--236.

\bibitem{ps1}
S. Ponnusamy and A. Sairam Kaliraj, Constants and characterization for certain classes of univalent harmonic mappings,
\textit{Mediterr. J. Math.} {\bf 12} (2015), 647--665.


\bibitem{sjr}
Y. Sun, Y.-P. Jiang and A. Rasila, On a subclass of close-to-convex harmonic mappings, \textit{Complex Var. Elliptic Equ.} {\bf 61} (2016), 1627--1643.

\bibitem{srj}
Y. Sun, A. Rasila and Y.-P. Jiang,
\newblock  Linear combinations of harmonic quasiconformal mappings convex in one direction,
\newblock \textit{Kodai Math. J.}  {\bf 39}  (2016), 366--377.

\bibitem{s2005}
T. J. Suffridge, Some special classes of conformal mappings, \textit{Handbook of complex analysis: geometric function theory}. Vol. 2, 309--338, Elsevier, Amsterdam, 2005.


\bibitem{wll}
Z.-G. Wang, Z.-H. Liu and Y.-C. Li, On the linear combinations of
harmonic univalent mappings, \textit{J. Math. Anal. Appl.} {\bf 400}
(2013), 452--459.

\bibitem{wsj1}
Z.-G. Wang, L. Shi and Y.-P. Jiang, On harmonic $K$-quasiconformal mappings associated with asymmetric vertical strips, \textit{Acta Math. Sin. (Engl. Ser.)} {\bf 31} (2015), 1970--1976.






\end{thebibliography}
\end{document}